%

\magnification=1200
\input amstex
\documentstyle{amsppt}

\overfullrule=0pt
\define\voj{Vojt\'a\v s}
\define\cc{\frak c}
\define\bb{\frak b}
\define\dd{\frak d}
\define\rr{\frak r}
\redefine\ss{\frak s}
\define\res{\restriction}
\define\m#1,#2{\text{Match}(#1,#2)}

\topmatter
\title
Reductions Between Cardinal Characteristics of the Continuum
\endtitle
\rightheadtext{Reductions Between Cardinal Characteristics}
\author
Andreas Blass
\endauthor
\address
Mathematics Dept., University of Michigan, Ann Arbor, 
MI 48109, U.S.A.
\endaddress
\email
ablass\@umich.edu
\endemail
\thanks Partially supported by NSF grant DMS-9204276 and NATO grant LG
921395.
\endgraf
Much of the material in this paper was presented at the set theory 
conference in Oberwolfach in October, 1993.
\endthanks
\abstract
We discuss two general aspects of the theory of cardinal
characteristics of the continuum, especially of proofs of inequalities
between such characteristics.  The first aspect is to express the
essential content of these proofs in a way that makes sense even in
models where the inequalities hold trivially (e.g., because the
continuum hypothesis holds).  For this purpose, we use a Borel version
of \voj's theory of generalized Galois-Tukey connections.  The second
aspect is to analyze a sequential structure often found in proofs of
inequalities relating one characteristic to the minimum (or maximum)
of two others.  \voj's max-min diagram, abstracted from such
situations, can be described in terms of a new, higher-type object in
the category of generalized Galois-Tukey connections.  It turns out to
occur also in other proofs of inequalities where no minimum (or
maximum) is mentioned.
\endabstract
\endtopmatter
\document

\head
1. Introduction
\endhead

Cardinal characteristics of the continuum are certain cardinal numbers
describing combinatorial, topological, or analytic properties of the
real line $\Bbb R$ and related spaces like ${}^\omega\omega$ and $\Cal
P(\omega)$.  Several examples are described below, and many more can
be found in \cite{4, 14}.  Most such characteristics, and all those under
consideration in this paper, lie between $\aleph_1$ and the
cardinality $\cc=2^{\aleph_0}$ of the continuum, inclusive.  So, if
the continuum hypothesis (CH) holds, they are equal to $\aleph_1$.
The theory of such characteristics is therefore of interest only when
CH fails.

That theory consists mainly of two sorts of results.  First, there are
equations and (non-strict) inequalities between pairs of
characteristics or sometimes between one characteristic and the
maximum or minimum of two others.  Second, there are independence
results showing that other equations and inequalities are not provable
in Zermelo-Fraenkel set theory (ZFC).  As examples of results of the
first sort, we mention in particular the work of Rothberger \cite{10} and
Bartoszy\'nski \cite1 relating the characteristics associated to
Lebesgue measure and Baire category; as examples of the second sort we
mention \cite2 and the earlier work cited there.  Many more examples
can be found in \cite{4, 14} and the references there.

A curious aspect of the proofs of inequalities (and of equations,
which we regard as pairs of inequalities) in this theory is that they
contain significant information whether or not CH holds, even though
CH makes the inequalities themselves trivial.  In other words, the
proofs establish additional information beyond the inequalities.
\voj\ \cite{15}  introduced a framework in which one can attempt to
formulate such additional information. He associates cardinal
characteristics to binary relations on the reals and shows that proofs
of inequalities between characteristics usually amount to the
construction of a suitable pair of functions between the domains and
the ranges of the corresponding relations.  He calls these pairs of
functions generalized Galois-Tukey connections, but for brevity we
shall call them morphisms.  The existence of such morphisms seems a
plausible candidate for the ``additional information'' established by
typical proofs of inequalities between cardinal characteristics.

It turns out, however, by a result of Yiparaki \cite{16}, that this
additional information is still trivial in the presence of CH.  We
shall show in Section~2 how to modify \voj's framework so as to
produce non-trivial results even in the presence of CH.  The key idea
is to add a definability (or absoluteness) requirement on the
functions that constitute a morphism.

Section~3 contains some information about Baire category that will be
used in subsequent examples.

Section~4 is concerned with a structure often found in proofs of
inequalities of the form $x\geq\min(y,z)$, a structure that \voj\
described with his max-min diagram.  We show that this diagram can be
neatly interpreted in terms of a construction, which we call {\sl
sequential composition\/}, on the relations associated to $y$ and $z$.
This construction allows us to analyze the ``flow of control'' in
proofs of such three-cardinal inequalities.

Section~5 is devoted to showing that sequential composition is
necessary for the results in Section~4.  In particular, certain
simpler compositions proposed by \voj\ are not adequate for these
results.

Finally, in Section~6, we discuss a situation where sequential
composition arises in the natural proof of an inequality involving
just two cardinal characteristics, not a maximum or minimum.

I thank Peter \voj\ and Janusz Pawlikowski for useful discussions of
the topic of this paper.

\head
2. Borel Galois-Tukey Connections
\endhead

To motivate \voj's framework, we define a few cardinal characteristics
of the continuum and discuss the proofs of some inequalities relating
them.  The characteristics to be used in this example are the {\sl
bounding number\/} $\bb$, the {\sl dominating number\/} $\dd$, the
{\sl unsplitting number\/} $\rr$, and the {\sl splitting number\/} $\ss$,
defined as follows.  (See \cite{4, 14} for more information about these
cardinals.)  For functions $f$ and $g$ from $\omega$ to $\omega$, we
write $f\leq^*g$ to mean that $f(n)\leq g(n)$ for all but finitely
many $n$.  A family $\Cal B\subseteq{}^\omega\omega$ of such functions
is {\sl unbounded\/} if there is no $g\in{}^\omega\omega$ such that
all elements of $\Cal B$ are $\leq^*g$.  The smallest possible
cardinality for an unbounded family is $\bb$.  A family $\Cal
D\subseteq{}^\omega\omega$ is {\sl dominating\/} if every
$g\in{}^\omega\omega$ is $\leq^*$ some $f\in\Cal D$.  The smallest
possible cardinality for a dominating family is $\dd$.  A subset $X$
of $\omega$ {\sl splits\/} another such set $Y$ if both $Y\cap X$ and
$Y-X$ are infinite.  A family $\Cal S\subseteq\Cal P(\omega)$ is a
{\sl splitting family\/} if every infinite $Y\subseteq\omega$ is split
by some member of $\Cal S$.  The smallest possible cardinality for a
splitting family is $\ss$.  A family $\Cal R$ of infinite subsets of 
$\omega$ is an {\sl unsplit family\/} (sometimes called {\sl
refining\/} or {\sl reaping\/}) if no single $X$ splits all the
members of $\Cal R$.  The smallest possible cardinality for an unsplit
family is $\rr$.

It is easy to see that all four of the cardinals just defined are
between $\aleph_1$ and $\cc$ inclusive and that $\bb\leq\dd$.  As an
example for future analysis, we give the proof of the less trivial yet
well-known inequality $\ss\leq\dd$.  Given a dominating family $\Cal
D$, we shall assign to each $f\in\Cal D$ a set $X_f\subseteq\omega$ in
such a way that all these $X_f$'s constitute a splitting family. This
will clearly suffice to prove $\ss\leq \dd$.  The construction is as
follows.  Given $f$, partition $\omega$ into finite intervals
$[a_0,a_1)$, $[a_1,a_2)$, \dots each satisfying $f(a_n)<a_{n+1}$.  To
be specific, let $a_0=0$ and $a_{n+1}=1+\max\{a_n,f(a_n)\}$.  Then let
$X_f$ be the union of the even-numbered intervals $[a_{2n},a_{2n+1})$.
To see that the $X_f$'s constitute a splitting family, let an
arbitrary infinite $Y\subseteq\omega$ be given, and let
$g:\omega\to\omega$ be the function sending each natural number $n$ to
the next larger element of $Y$.  As $\Cal D$ is dominating, it
contains an $f\geq^*g$.  In the construction of $X_f$ we have, for all
sufficiently large $n$, that the next element of $Y$ after $a_n$ is
$g(a_n)\leq f(a_n)<a_{n+1}$ and therefore lies in the interval
$[a_n,a_{n+1})$.  So $Y$ meets all but finitely many of these
intervals and therefore contains infinitely many members of $X_f$ and
infinitely many members of $\omega-X_f$.  That is, $X_f$ splits $Y$,
as required.

We now present \voj's framework for describing cardinal
characteristics and proofs of inequalities, using the preceding proof
as an example.  First, each characteristic was defined as the smallest
possible cardinality for a set $\Cal Z$ of reals such that every real
is related in a certain way to one in $\Cal Z$.  More precisely, in
each case we had a triple $\bold A=(A_-,A_+,A)$ of two sets $A_\pm$
and a binary relation $A$ between them such that the characteristic is
$$
\Vert(A_-,A_+,A)\Vert=\min\{|\Cal Z|:\,\Cal Z\subseteq A_+\text{ and
}
\forall x\in A_-\,\exists z\in\Cal Z\, A(x,z)\},
$$
which we call the {\sl norm\/} of $\bold A$.  For example, $\dd$ is
the norm of $({}^\omega\omega,{}^\omega\omega,\leq^*)$, and $\ss$ is
the norm of $(\Cal P_\infty(\omega),\Cal P(\omega),\,\text{is split
by})$.  (Here $\Cal P_\infty$ means the family of infinite subsets.)
To decribe $\bb$ and $\rr$ as norms, it suffices to take the
descriptions for $\dd$ and $\ss$ and {\sl dualize\/} them in the
following sense: interchange $A_-$ with $A_+$, and replace $A$ with
the complement of the converse relation.  In general, we write
$$
(A_-,\hskip0.01cm A_+,A)^\bot=(A_+,A_-,\{(z,x):\,\text{not }A(x,z)\}.
$$
(\voj\ calls the norms of $\bold A$ and of its dual the dominating
and bounding numbers of $\bold A$, respectively, by analogy with
the example of $\dd$ and $\bb$ above.)  We sometimes call triples
$(A_-,A_+,\hskip0.01cm A)$ relations, although strictly speaking it is only the
third component $A$ that is a relation.  To avoid trivialities, we
shall tacitly assume that our relations have $A_\pm\neq\emptyset$,
that each element of $A_-$ is $A$-related to some element of $A_+$ (so
$\Vert\bold A\Vert$ is defined), and that not all elements of $A_-$
are $A$-related to any single element of $A_+$ (so $\Vert\bold
A^\bot\Vert$ is defined).  These assumptions amount to requiring the
norms of both $\bold A$ and $\bold A^\bot$ to be at least 2.

The proof of $\ss\leq\dd$ presented above consists of three parts.
First, there was a construction $\xi_+:{}^\omega\omega\to\Cal
P(\omega)$ sending each $f$ to $X_f$.  Second, there was a (simpler)
construction $\xi_-:\Cal P_\infty(\omega)\to{}^\omega\omega$ sending
each $Y$ to the function $\xi_-(Y)=g$ that maps $n$ to the next larger
element of $Y$.  Finally, there was the verification that if
$g\leq^*f$ then $Y$ is split by $X_f$.  Abstracting this structure, we
obtain \voj's notion of a generalized Galois-Tukey connection; we call
it simply a morphism and write it in the opposite direction to \voj's.

\definition{Definition}
A {\sl morphism\/} $\xi$ from $(A_-,A_+,A)$ to $(B_-,B_+,B)$ is a pair of
functions $\xi_-:B_-\to A_-$ and $\xi_+:A_+\to B_+$ such that, for all
$b\in B_-$ and all $a\in A_+$,
$$
A(\xi_-(b),a)\implies B(b,\xi_+(a)).
$$
\enddefinition

Our convention for the direction of morphisms was chosen partly to
work well with other uses of the same category \cite{8, 3} and partly so
that the direction of a morphism agrees with the direction of the
implication displayed in the definition.

The existence of a morphism $\xi$ from $(A_-,A_+,A)$ to $(B_-,B_+,B)$
immediately implies the norm inequality $\Vert(A_-,A_+,A)\Vert\geq
\Vert(B_-,B_+,B)\Vert$.  Indeed, if $\Cal Z$ is as in the definition
of norm for $(A_-,A_+,A)$, then $\xi_+(\Cal Z)$ has no greater
cardinality and serves the same purpose in the definition of the norm
of $(B_-,B_+,B)$.  Because of this, we write
$\bold A\geq\bold B$ to indicate the existence of such a
morphism.  (Our inequalities, unlike our morphisms, go in the same
direction as \voj's.)

A morphism from $\bold A$ to $\bold B$ becomes, just by interchanging
its two components, a morphism in the opposite direction between the
dual objects.  Thus, for example, the proof above of $\ss\leq\dd$,
exhibiting a morphism from $({}^\omega\omega,{}^\omega\omega,\leq^*)$
to $(\Cal P_\infty(\omega),\Cal P(\omega),\,\text{is split by})$, also
exhibits a morphism from $(\Cal P(\omega),\Cal
P_\infty(\omega),\,\text{does }\allowmathbreak\text{not
}\allowmathbreak\text{split})$ to
$({}^\omega\omega,{}^\omega\omega,\ngeq^*)$ and thus proves that
$\bb\leq\rr$.  (In this form, the inequality is essentially due to
Solomon \cite{11}.)

One can similarly exhibit morphisms that capture the combinatorial
content of the proofs of a great many other inequalities between
cardinal characteristics.  This applies both to trivial inequalities
like $\bb\leq\dd$ and deep results like Bartoszy\'nski's theorem that
the additivity of category is at least the additivity of measure.  See
\cite5 for a presentation of Bartoszy\'nski's theorem that makes the
morphism explicit, and see \cite{15} for more examples.

For future reference, we mention that the converse, $\dd\leq\ss$, of
the inequality proved above is known not to be provable in ZFC (though
of course it holds in some models of ZFC, for example models of CH).
It fails, for example, in the model obtained from a model of CH by
adding a set $C\subseteq{}^\omega\omega$ of $\aleph_2$ Cohen reals.
In this model, $\dd=\aleph_2$ because any $\aleph_1$ reals lie in a
submodel generated by $\aleph_1$ members of $C$ and therefore cannot
dominate the other members of $C$.  On the other hand, $\ss=\aleph_1$
because the set of ground model reals is non-meager (i.e., of second
Baire category) in the extension and any non-meager subset of $\Cal
P(\omega)$ is a splitting family (because the reals that fail to split
any particular $Y$ form a meager set).

We shall call an equation or (non-strict) inequality between cardinal
characteristics ``correct'' if it is provable in ZFC and ``incorrect''
if it is independent of ZFC.  (It cannot be refutable in ZFC, because
it holds in models of CH; recall that we deal only with
characteristics that lie between $\aleph_1$ and $\cc$.)  In models of
CH, the difference between correct and incorrect inequalities is
hidden, because all the inequalities are true there.  Nevertheless,
it seems reasonable to say, even in such models, that correct
inequalities, like $\ss\leq\dd$, hold for understandable,
combinatorial reasons while incorrect inequalities, like $\dd\leq\ss$,
hold only because CH ``happens'' to be true.  Can one make
mathematical sense of such statements?

\voj's theory, in particular the fact that proofs of inequalities
between characteristics usually exhibit morphisms between the
corresponding relations, suggests an affirmative answer to this
question.  The understandable, combinatorial reason for a correct
inequality is given by the morphism.  So one might hope that there
are, even in the presence of CH, no morphisms corresponding to
incorrect inequalities. Then, when incorrect inequalities hold in a
model, this would not be because of good reasons (i.e., morphisms) but
because of other properties of the model (like CH).

A theorem of Yiparaki \cite{16}, Chapter~5, dashes this hope.  She
shows that, if
$$
\Vert\bold A\Vert=|A_+|=\Vert\bold B^\bot\Vert=|B_-|,
$$ 
then there is a morphism from $\bold A$ to $\bold B$.  In
particular, in all our examples there is such a morphism if CH holds,
because all the cardinals in the displayed equation are then equal to
$\aleph_1$.  So models of CH not only satisfy all inequalities,
correct or incorrect, they also contain morphisms to justify all these
inequalities.

These morphisms, however, are highly non-constructive; their
definition involves a multitude of arbitrary choices.  We
propose therefore to eliminate them by working not with arbitrary
morphisms but with well-behaved ones.  To be specific, we restrict our
attention to objects $(A_-,A_+,A)$ in which $A_\pm$ are Borel sets of
reals and $A$ is a Borel relation between them, and we consider only
morphisms $\xi$ both of whose components $\xi_\pm$ are Borel
functions.  The objects and morphisms considered above, in connection
with the definitions of $\dd$ and $\ss$ (and their duals) and the
proof of $\ss\leq\dd$ (and its dual) are all Borel in this sense.
Yiparaki's proof, on the other hand, involves non-Borel morphisms.  

As an example of what the restriction to Borel morphisms accomplishes,
we show that there is no Borel morphism corresponding to the incorrect
inequality $\dd\leq\ss$.

\proclaim{Proposition 1}
There is no Borel morphism
$$
(\xi_-,\xi_+):(\Cal P_\infty(\omega),\Cal P(\omega),
\text{is split by})
\to({}^\omega\omega,{}^\omega\omega,\leq^*).
$$
\endproclaim

\demo{Proof}
Suppose we had such a morphism $(\xi_-,\xi_+)$ consisting of two
Borel maps.  Let $c\in{}^\omega\omega$ be Cohen-generic over the
universe $V$ (in some Boolean extension), and by abuse of notation
write $(\xi_-,\xi_+)$ also for the pair of Borel maps in $V[c]$ having
the same codes as the original $(\xi_-,\xi_+)$ had in $V$.  Since the
ground model reals remain a non-meager and hence splitting set in the
extension, $\xi_-(c)$ is split by some $X$ in the ground
model. Because $\xi$ is a morphism, it follows that $c\leq^*\xi_+(X)$.
But the ground model contains $X$ and the code for $\xi_+$ and
therefore also $\xi_+(X)$.  This is absurd, as no real from the ground
model can dominate a Cohen real.
\qed\enddemo

More generally, we can show that there are no Borel morphisms between
Borel relations when the corresponding inequality of characteristics
can be violated by forcing.  To express this precisely, suppose we
have Borel $\bold A$ and $\bold B$ such that some notion of forcing
$P$ forces $\Vert\bold A\Vert<\Vert\bold B\Vert$, where we have, as
above, abused notation by writing $\bold A$ and $\bold B$ for the
objects in the forcing extension having the same Borel codes as the
original $\bold A$ and $\bold B$ in the ground model.  Then there is
no Borel morphism (in the ground model) $\xi:\bold A\to\bold B$.
Indeed, for $\xi$ to be such a morphism would be a $\Pi^1_1$ assertion
about the Borel codes of $\xi$, $\bold A$, and $\bold B$.  Such
assertions are preserved by forcing extensions.  But $P$ forces that
there is no such morphism, because it forces the opposite ordering of
the norms.  (Note that we proved Proposition~1 by adding a single
Cohen real, whereas the general argument just given would involve
adding at least $\aleph_2$ Cohen reals.)

We emphasize that, in the situation of the preceding paragraph, the
inequality $\Vert\bold A\Vert\geq\Vert\bold B\Vert$ may well be true
in the ground model, but there cannot be a Borel morphism causing it.

It should be noted that independence proofs for inequalities between
cardinal characteristics (e.g., \cite{2, 14}) typically take the form
considered above, i.e., they produce a $P$ forcing a strict inequality
in the opposite direction.  In many cases, the proofs in the
literature involve not only a forcing construction but some
requirements on the ground model over which the forcing is done.  But
these requirements, usually CH or Martin's axiom, can themselves be
forced, so there is a $P$ as in the discussion above.  

Summarizing, we have that incorrect inequalities between cardinal
characteristics, though they may be true in some models and may be
given by morphisms in some models, are not given by Borel morphisms in
any model, provided their incorrectness can be established by a
forcing argument.

\head
3. Baire Category
\endhead

In Section~4, we shall use some cardinal characteristics related to
Baire category as an example to motivate and illustrate the operation
of sequential composition of relations.  In preparation for this, we
devote the present section to introducing the notation and preliminary
results needed to make the later discussion proceed smoothly.  Along
the way, we shall see a few more examples of the structures discussed
in Section~2.

In what follows, we shall use the Cantor space ${}^\omega2$, with its
usual (product) topology, as the underlying space in all our
discussions of Baire category.  One of the characteristics that we
shall need is the additivity of Baire category, ${\bold{add}}(B)$, the
smallest number of meager (= first category) sets whose union is not
meager.  Like the characteristics discussed in the preceding section,
${\bold{add}}(B)$ lies between $\aleph_1$ and $\cc$ inclusive.  It can
be described as the norm of $(B,B,\nsupseteq)$, where $B$ is the
collection of meager subsets of ${}^\omega2$.  This description is not
amenable, as it stands, to the Borel considerations of the preceding
section, because $B$ is not a set of reals but a set of sets of reals.
One can, however, replace $B$ with the collection of meager $F_\sigma$
sets, i.e., countable unions of nowhere-dense closed sets; this does
not affect ${\bold{add}}(B)$, since every meager set is included in a
meager $F_\sigma$ set.  It is easy to code meager $F_\sigma$ sets by
reals in such a way that the set of codes is Borel, and with some
additional work one can arrange that the relations we need, like
$\nsupseteq$, are also Borel in the codes.  We omit the details of
this since we shall soon introduce a different way of viewing
${\bold{add}}(B)$ for which these matters are easier to handle.

The other Baire category characteristic that we shall need is the
covering number, ${\bold{cov}}(B)$, the smallest number of meager sets
needed to cover ${}^\omega2$.  This is the norm of $({}^\omega2,
B,\in)$, and again we can replace $B$ by the collection of meager
$F_\sigma$ sets or by the collection of their codes.  It is obvious
that ${\bold{add}}(B)\leq{\bold{cov}}(B)$.  As expected, this trivial
inequality corresponds to a trivial morphism, with
$\xi_-:B\to{}^\omega2$ sending any meager set to some real not in it
and with $\xi_+:B\to B$ being the identity map.  (The map $\xi_-$
involves an arbitrary choice, but from a code for a meager $F_\sigma$
set one can obtain, in a Borel fashion, a specific real not in that
set, just by following the proof of the Baire category theorem.)

It will be convenient to use a particular, easily coded basis for the
ideal of meager sets.  To introduce this basis, we first define a {\sl
chopped real\/} to be a pair $(x,\Pi)$, where $x\in{}^\omega2$ and
where $\Pi$ is a partition of $\omega$ into finite intervals
$I_0=[a_0,a_1)$, $I_1=[a_1,a_2)$, \dots, with $0=a_0<a_1<a_2<\dots$.
The idea is that a real $x:\omega\to2$ has been chopped into finite
pieces by the partition $\Pi$ of its domain.  We say that a real
$y\in{}^\omega2$ {\sl matches\/} the chopped real $(x,\Pi)$ if $y$
agrees with $x$ on infinitely many of the intervals of $\Pi$, i.e., if
there are infinitely many $n\in\omega$ such that $y\res I_n=x\res
I_n$.  We write $\m x,\Pi$ for the set of all $y$ that match
$(x,\Pi)$.  It is easy to verify that $\m x,\Pi$ is a dense $G_\delta$
set.  Talagrand \cite{12}, in obtaining a combinatorial description of
meager filters, proved that every dense $G_\delta$ set in ${}^\omega2$
has a subset of the form $\m x,\Pi$; equivalently, every meager set is
included in the complement of $\m x,\Pi$ for some chopped real.

Thus, we can use chopped reals and the corresponding matching sets
instead of arbitrary meager sets or meager $F_\sigma$ sets in
describing the cardinal characteristics of Baire category.
Specifically, let us define $CR$ to be the set of chopped reals and
define 
$$
\bold U=(CR,{}^\omega2,\text{matches}^\smile),
$$ 
where ${}^\smile$ means to take the converse of a relation.  Then
the norm of the dual of $\bold U$ is the smallest size for a family of
chopped reals such that no single real matches them all.
Equivalently, by Talagrand's result, it is the smallest number of
dense $G_\delta$ sets with empty intersection.  Taking the complements
of those sets, we find that this is precisely the covering number.  So
we have ${\bold{cov}}(B)=\Vert\bold U^\bot\Vert$.  (The norm of $\bold
U$ itself is the smallest cardinality of a non-meager set of reals,
called the {\sl uniformity\/} of category; it is important in its own
right, but we shall not need it here.  The reason for defining $\bold
U$ as we did, rather than dually, is to avoid excessive negations and
to avoid some dualizations in the next section.)

For a similar description of the additivity of category, we need to
describe, in terms of chopped reals, the inclusion relation between
the corresponding dense $G_\delta$ sets.  We say that one chopped real
$(x,\Pi)$ {\sl engulfs\/} another $(x',\Pi')$ if all but finitely many
intervals of $\Pi$ include intervals of $\Pi'$ on which $x$ and $x'$
agree.  

\proclaim{Lemma}
$\m x,\Pi\subseteq\m x',\Pi'$ if and only if $(x,\Pi)$ engulfs
$(x',\Pi')$. 
\endproclaim

\demo{Proof}
First, suppose $(x,\Pi)$ engulfs $(x',\Pi')$ and $y$ matches
$(x,\Pi)$.  So there are infinitely many intervals $I$ of $\Pi$ on
which $y$ agrees with $x$.  Except for finitely many, each such $I$
includes an interval $J$ of $\Pi'$ on which $x$ and $x'$ agree.  Thus
we get infinitely many intervals $J$ of $\Pi'$ on which $y$ and $x'$
agree, so $y$ matches $(x',\Pi')$.

Conversely, suppose $(x,\Pi)$ does not engulf $(x',\Pi')$, so there
are infinitely many intervals $I$ of $\Pi$ that contain no interval of
$\Pi'$ on which $x$ and $x'$ agree.  Discarding some of these
intervals $I$, we can arrange that no interval of $\Pi'$ meets more
than one of them.  Define $y\in{}^\omega2$ by making it agree with $x$
on the union of these intervals $I$ and making it disagree with $x'$
everywhere else.  Then $y$ matches $(x,\Pi)$ but not $(x',\Pi')$.
\qed\enddemo

Define
$$
\bold V=(CR,CR,\text{is engulfed by}).
$$
Then the norm of the dual of $V$ is the minimum number of chopped
reals such that no single real engulfs them all.  Equivalently, by
Talagrand's result and the lemma, it is the minimum number of dense
$G_\delta$ sets such that no single dense $G_\delta$ set is included
in them all.  Taking complements, we find that this is just the
additivity of category.  So ${\bold{add}}(B)=\Vert\bold V^\bot\Vert$.
(The norm of $\bold V$ itself is the {\sl cofinality\/} characteristic
of Baire category.)

The trivial inequality ${\bold{add}}(B)\leq{\bold{cov}}(B)$ is given
by a trivial morphism $\bold U^\bot\to\bold V^\bot$ or equivalently
$\bold V\to\bold U$ in which one component is the identity map $CR\to
CR$ while the other component $CR\to{}^\omega2$ sends a chopped real
to its first component (forgetting the partition and keeping only the
real).  

A non-trivial inequality \cite7, namely ${\bold{add}}(B)\leq\bb$, is
also easy to see from this point of view.  Notice that $\bb$ is, as
discussed in Section~2, the norm of $\bold W^\bot$, where
$$
\bold W=({}^\omega\omega,{}^\omega\omega,\leq^*).
$$
So to prove the inequality in question, it suffices to exhibit a
morphism $\bold W^\bot\to\bold V^\bot$ or equivalently $\xi:\bold
V\to \bold W$.  Define $\xi_-$ to map any $f\in{}^\omega\omega$ to the
chopped real $(0,\Pi)$ where 0 is the identically zero function
$\omega\to2$ and where $\Pi$ is chosen so that, for any $n\in\omega$,
$f(n)$ is no more than one interval past $n$ (i.e., if $n\in I_k$ then
$f(n)\in I_l$ for some $l\leq k+1$).  (To get a Borel map $\xi_-$, the
intervals should be chosen in a canonical manner, for example by
choosing each endpoint $a_n$ as small as possible subject to the
constraints that $f(n)$ be at most one interval past $n$ and that all
the intervals be nonempty.)  Define $\xi_+$ to map any chopped real
$(x,\Pi)$ to the function in ${}^\omega\omega$ sending each
$n\in\omega$ to the right endpoint of the interval after the one that
contains $n$.  It is straightforward to verify that, if $\xi_-(f)$ is
engulfed by $(x,\Pi)$, then $f\leq^*\xi_+(x,\Pi)$, so $\xi$ serves as
the required morphism.

The inequalities above combine to give
${\bold{add}}(B)\leq\min\{{\bold{cov}}(B),\bb\}$.  In fact, equality
holds here \cite{13, 7}, but the converse inequality cannot be
separated into two simpler inequalities each proved by exhibiting a
morphism.  On the contrary, this converse inequality, like any
inequality of the form $x\geq\min\{y,z\}$ (or $x\leq\max\{y,z\}$)
involves an interaction of all three cardinals.  In the next section,
we shall review the proof of
${\bold{add}}(B)\geq\min\{{\bold{cov}}(B),\bb\}$ and use it to
motivate a construction that allows such proofs to be presented as
morphisms between suitable relations.

\head
4. Sequential Composition
\endhead

To treat inequalities involving the maximum or minimum of two cardinal
characteristics, it is natural to seek a construction which, given two
relations $\bold A$ and $\bold B$, produces another relation $\bold C$
with $\Vert C\Vert=\max\{\Vert A\Vert,\Vert B\Vert\}$ and $\Vert
C^\bot\Vert=\min\{\Vert A^\bot\Vert,\allowmathbreak\Vert
B^\bot\Vert\}$.  Such a construction is the product $\bold
A\times\bold B$ given by letting $C_-$ be the disjoint union
$A_-\sqcup B_-$, letting $C_+$ be the product $A_+\times B_+$, and
defining $C(x,(a,b))$ to hold when either $x\in A_-$ and $A(x,a)$ or
$x\in B_-$ and $B(x,b)$.  This is the product in the
category-theoretic sense with respect to our definition of morphisms.
(With \voj's convention it is, as he points out, the coproduct.)  It
would be pleasant if proofs of three-cardinal inequalities, like the
${\bold{add}}(B)\geq\min\{{\bold{cov}}(B),\bb\}$ mentioned at the end
of the preceding section, could be presented as constructions of
morphisms from a product, and \voj\ asks (\cite{15}) whether this can be
done, after noting that the usual proof has a more complicated
structure (described below).

An earlier (preprint) version of \cite{15} contained a different sort of
product, which we shall call the {\sl old product\/} to distinguish it
from the (categorical) product in the preceding paragraph.  The old
product $\bold C$ of $\bold A$ and $\bold B$ has $C_-=A_-\times B_-$,
$C_+=A_+\times B_+$, and $C((x,y),(a,b))$ if and only if both $A(x,a)$
and $B(y,b)$.  The norm of $\bold C$ and its dual are the maximum and
minimum of the norms of the factors and their duals, respectively, as
long as the norms are infinite.  (More precisely, $\Vert\bold
C^\bot\Vert=\min\{\Vert\bold A^\bot\Vert,\Vert\bold B^\bot\Vert\}$ and
$\max\{\Vert\bold A\Vert,\Vert\bold B\Vert\}\leq\Vert\bold
C\Vert\leq\Vert\bold A\Vert\cdot\Vert\bold B\Vert$.  The last two
inequalities can both be strict, but of course only when the norms
involved are finite.  For example, $(3,3,\neq)$ has norm 2 but its old
product with itself has norm 3.)  In the preprint version of \cite{15},
\voj\ asked whether proofs of certain three-cardinal inequalities
could be presented as constructions of morphisms from an old product.

We shall show that the answer to his question is negative for both the
product and the old product if we require morphisms to consist of Borel
maps.  (Without this requirement, Yiparaki's result applies and shows
that an affirmative answer is consistent, being true in models of CH.)
Thus the more complicated structure, described by \voj\ in his max-min
diagram, is essential.  But we shall see that this more complicated
structure can also be described in terms of a construction of a
suitable (more complicated) $\bold C$ from $\bold A$ and $\bold B$.

To motivate this construction and set the stage for the proof of its
necessity, we review the proof that
${\bold{add}}(B)\geq\min\{{\bold{cov}}(B),\bb\}$, using the notation
and machinery of the preceding section.

Let $(x_\alpha, \Pi_\alpha)$ for $\alpha<\kappa$ be
$\kappa<\min\{{\bold{cov}}(B),\bb\}$ chopped reals.  We must produce a
single chopped real $(y,\Theta)$ engulfing them all.  Since
$\kappa<{\bold{cov}}(B)$, fix a real $y$ matching all the $(x_\alpha,
\Pi_\alpha)$.  This will be the first component of the chopped real we
seek; it remains to produce $\Theta$ such that, for each
$\alpha<\kappa$, for all but finitely many blocks $J$ of $\Theta$,
there is a block $I$ of $\Pi_\alpha$ such that $I\subseteq J$ and
$x_\alpha\res I=y\res I$.

For each $\alpha$, define $f_\alpha:\omega\to\omega$ by letting
$f_\alpha(n)$ be the right endpoint of the next interval of
$\Pi_\alpha$, after the one containing $n$, on which $y$ agrees with
$x_\alpha$.  (Such an interval exists because $y$ matches $(x_\alpha,
\Pi_\alpha)$.)  As $\kappa<\bb$, fix some $g:\omega\to\omega$
eventually majorizing every $f_\alpha$.  Then choose $\Theta$ so that,
for each of its intervals, $g$ of the left endpoint is smaller than
the right endpoint.  For each $\alpha<\kappa$, for all but finitely
many intervals $J=[a,b]$ of $\Theta$, we have $f_\alpha(a)\leq g(a)<b$
and therefore, by definition of $f_\alpha$,  $y$ agrees with
$x_\alpha$ on some interval of $\Pi_\alpha$ that starts after $a$ and
ends before $b$ and thus is contained in $J$.  This shows that
$(y,\Theta)$ is as required.

Let us describe this proof in terms of the relations $\bold U$, $\bold
V$, and $\bold W$ from the preceding section.  Recall that 
$$\align
\bold U&=(CR,{}^\omega2,\text{matches}^\smile),\\
\bold V&=(CR,CR,\text{is engulfed by}),\\
\bold W&=({}^\omega\omega,{}^\omega\omega,\leq^*),
\endalign$$
so $\Vert\bold U^\bot\Vert={\bold{cov}}(B)$, $\Vert\bold V^\bot\Vert=
{\bold{add}}(B)$, and $\vert\bold W^\bot\Vert=\bb$.  (Thus the
inequality ${\bold{add}}(B)\geq\min\{{\bold{cov}}(B),\bb\}$ should
correspond to a morphism from some sort of ``product'' of $\bold U$
and $\bold W$ to $\bold V$; the problem is to define an appropriate
sort of product.)

The proof of ${\bold{add}}(B)\geq\min\{{\bold{cov}}(B),\bb\}$ above
began by regarding the given elements $(x_\alpha,\Pi_\alpha)$ of $V_-$
(for which we needed to find a $V$-related element of $V_+$) as
elements of $U_-$ and finding a $y\in U_+$ that is $U$-related to them
all.  Thus, the proof implicitly used the identity map to convert
elements of $V_-$ into elements of $U_-$.  (In other applications, a
non-trivial map will occur here.)  The next step was to
define, from $y$ and $(x_\alpha,\Pi_\alpha)$, the element $f_\alpha\in
W_-$; so the proof uses a map $V_-\times U_+\to W_-$.  (Our definition
of $f_\alpha$ presupposed that $y$ matches $(x_\alpha,\Pi_\alpha)$; to
get a total map, let $f_\alpha$ be identically zero if $y$ fails to
match $(x_\alpha,\Pi_\alpha)$.)  Finally, from
$y$ and an element $g\in W_+$ that is $W$-related to each $f_\alpha$,
we produced the required $(y,\Theta)\in V_+$ that is $V$-related to
each $(x_\alpha,\Pi_\alpha)$.  The constructions in the proof can be
summarized as three maps (where we have omitted the subscripts
$\alpha$ from $x$, $\Pi$, and $f$)
$$\align
\alpha:V_-\to U_-&:(x,\Pi)\mapsto(x,\Pi),\\
\beta:V_-\times U_+\to W_-&:((x,\Pi),y)\mapsto f,
\text{ and}\\
\gamma:U_+\times W_+\to V_+&:(y,g)\mapsto(y,\Theta).
\endalign$$
Their key property is that from $U(\alpha(x,\Pi),y)$ and
$W(\beta((x,\Pi),y),g)$ we were able to infer
$V((x,\Pi),\gamma(y,g))$.

Notice that if $\beta$ were a function of only $(x,\Pi)$ rather than
both $(x,\Pi)$ and $y$, then the situation above would precisely
describe a morphism $\xi$ from the old product of $\bold U$ and $\bold
W$ to $\bold V$.  Indeed, $\alpha$ and $\beta$ could then serve as the
two components of $\xi_-:V_-\to U_-\times W_-$ and $\gamma$ could
serve as $\xi_+$.

But $\beta$ in the proof definitely needs $y$ as an argument, so this
proof is not described by a morphism from an old product.  It is even
less describable in terms of the (new, categorical) product; that
would require each element of $V_-$ to have an image in $U_-$ {\sl
or\/} in $W_-$, not both.

The awkward $U_+$ in the domain of $\beta$ can be moved to the
codomain by considering what category theorists call the exponential
adjoint and computer scientists call currying:
$$
\hat\beta:V_-\to{}^{U_+}W_-:(x,\Pi)\mapsto(y\mapsto\beta((x,\Pi),y)).
$$
The key property of $\alpha$, $\beta$, and $\gamma$ can, of course, be
trivially rewritten in terms of $\alpha$, $\hat\beta$, and $\gamma$.
The result is that the first two of these are the components of
$\xi_-$ and the third is $\xi_+$ for a morphism $\xi$ to $\bold V$
from the following compound of $\bold U$ and $\bold W$, which we call
their {\sl sequential composition\/} (cf.\ also \cite3):
$$
\bold U;\bold W=(U_-\times{}^{U_+}W_-, U_+\times W_+,
\{((a,\rho),(u,w))\mid U(a,u)\text{ and }W(\rho(u),w)\}.
$$

Notice that, quite generally, a morphism from a sequential composition
to another relation can be regarded (via exponential adjointness) as
consisting of three functions that enjoy the key property of $\alpha$,
$\beta$, and $\gamma$ discussed above.  Such a triple of functions is
precisely what \voj\ \cite{15} describes with his max-min diagram.

The following proposition records the connection between sequential
composition and maxima and minima of cardinal characteristics.  Note
that this proposition has a product of cardinals where the product of
relations has a maximum of cardinals and the old product has (as
discussed above) something between the maximum and the product of
cardinals.  Of course, for infinite cardinals (the case we're
interested in) all these cardinals coincide.

\proclaim{Proposition 2}
For any relations $\bold A$ and $\bold B$, we have 
$$
\Vert\bold A;\bold B\Vert=\Vert\bold A\Vert\cdot\Vert\bold B\Vert
\quad\text{and}\quad
\Vert(\bold A;\bold B)^\bot\Vert=
\min\{\Vert\bold A^\bot\Vert,\Vert\bold B^\bot\Vert\}.
$$
\endproclaim

\demo{Proof}
We prove only the least trivial part, namely the $\geq$ half of the
first equation.  By definition, $\Vert\bold A;\bold B\Vert$ is the
smallest possible cardinality for a set $S\subseteq A_+\times B_+$
such that, for every element $x\in A_-$ and every function
$\rho:A_+\to B_-$, some $(a,b)\in S$ satisfies $A(x,a)$ and
$B(\rho(a),b)$.  Fix such an $S$ of minimum cardinality; we must show
that this cardinality is at least the product of the norms of $\bold
A$ and $\bold B$.

For each $a\in A_+$, let $S_a=\{b\in B_+\mid(a,b)\in S\}$.  Let
$X=\{a\in A_+\mid |S_a|\geq\Vert\bold B\Vert\}$.  It will suffice to
prove that $X$ has cardinality at least $\Vert\bold A\Vert$, for then
we have at least $\Vert\bold A\Vert$ elements $a$ each $S$-related to
at least $\Vert\bold B\Vert$ elements $b$ and so we have at least
$\Vert\bold A\Vert\cdot\Vert\bold B\Vert$ pairs $(a,b)\in S$.

So suppose toward a contradiction that $|X|<\Vert\bold A\Vert$.  By
definition of norm, we can fix some $x\in A_-$ that is not $A$-related
to any element of $X$.  Using the definition of norm again along with
the definition of $X$, we can choose, for each $a\in A_+-X$, an
element $\rho(a)\in B_-$ that is not $B$-related to any $b\in S_a$.
Extend $\rho$ arbitrarily to a function $\rho:A_+\to B_-$.  By our
choice of $S$, there is some $(a,b)\in S$ with $A(x,a)$ and
$B(\rho(a),b)$.  From $A(x,a)$ we infer, by our choice of $x$, that
$a\notin X$.  But then from $B(\rho(a),b)$ and our choice of $\rho$ it
follows that $b\notin S_a$.  That contradicts $(a,b)\in S$.
\qed\enddemo

By virtue of this proposition and the general properties of morphisms,
we see that the existence of a morphism $\bold A;\bold B\to\bold C$
implies $\Vert\bold C\Vert\leq\Vert\bold A\Vert\cdot\Vert\bold
B\Vert\,(=\max\{\Vert\bold A\Vert,\Vert\bold B\Vert\}$ if the norms
are infinite) and $\Vert\bold C^\bot\Vert\geq\min\{\Vert\bold
A^\bot\Vert,\Vert\bold B^\bot\Vert\}$.  In other words, morphisms from
sequential compositions provide a way to present proofs of
three-cardinal inequalities relating (in the non-trivial direstion)
one cardinal to the maximum or minimum of two others.  The example of
${\bold{add}}(B)\geq\min\{{\bold{cov}}(B),\bb\}$ shows that this sort
of presentation occurs naturally in one example; other examples are
given in \voj's discussion of the max-min diagram \cite{15}, and another
example (involving a sequential composition one of whose factors is an
old product) is discussed in \cite3.

One unpleasant aspect of sequential composition needs to be discussed.
Because of the function-space construction in the definition of
sequential composition, even if $\bold A$ and $\bold B$ are, as we
advocated in Section~2, Borel relations between Borel sets of reals,
$\bold A;\bold B$ will not be of this form, simply because its domain
is not a set of reals at all but rather a set one type higher.  So it
does not yet make sense to talk about Borel morphisms involving
sequential compositions.

Fortunately, when we deal (as we did above) with morphisms $\xi$ {\sl
from\/} a sequential composition $\bold A;\bold B$ to some $\bold C$,
then the troublesome (higher type) part is $\xi_-:C_-\to
A_-\times{}^{A_+}B_-$ or more precisely its second component
$C_-\to{}^{A_+}B_-$.  We can declare such a function to be Borel if
and only if its exponential adjoint $C_-\times A_+\to B_-$ is Borel.
In this sense, the morphism involved in our proof of
${\bold{add}}(B)\geq\min\{{\bold{cov}}(B),\bb\}$ is Borel, and so are
the morphisms involved in the other proofs to which we alluded two
paragraphs ago.

The exponential adjoint trick does not work for morphisms {\sl to\/} a
sequential composition, for then the function space ${}^{A_+}B_-$
occurs in the domain rather than the codomain of a function.  The only
sensible way to talk about such functions being Borel seems to be to
work not with the space of all functions $A_+\to B_-$ but with the
subspace of Borel functions, or rather with the set of Borel codes for
such functions.  We postpone further discussion of this until after we
see, in the next section, an occurrence ``in nature'' of something
that ought to be a Borel morphism to a sequential composition.

\head
5. Non-existence of Borel morphisms
\endhead

The purpose of this section is to show that, in the proof of the
inequality ${\bold{add}}(B)\geq\min\{{\bold{cov}}(B),\bb\}$, one
cannot replace the sequential composition $\bold U;\bold W$ with the
product or the old product.  In fact, we shall see that $\bold V$ is
in some sense equivalent to $\bold U;\bold W$.  We begin with a
proposition showing that products (old or new) do not suffice.  It
also shows that the order of the components in a sequential
composition is essential.

\proclaim{Proposition 3}
Let $\bold U$, $\bold V$, and $\bold W$ be as in the preceding two
sections.  There is no Borel morphism to $\bold V$ from the product of
$\bold U$ and $\bold W$, nor from their old product, nor from $\bold
W;\bold U$.
\endproclaim

\demo{Proof}
Notice first that there are morphisms from $\bold W;\bold U$ to the
old product and from the old product to the product.  Both morphisms
are the identity on the $+$-components $W_+\times U_+$.  On the
$-$-components, we have $W_-\sqcup U_-\to W_-\times U_-\to
W_-\times{}^{W_+}U_-$, where the first map is defined by fixing
elements $w_0\in W_-$ and $u_0\in U_-$ and sending any $w\in W_-$ to
$(w,u_0)$ and any $u\in U_-$ to $(w_0,u)$, and the second map is the
identity on the first component $W_-$ while on the second component it
sends any $u\in U_-$ to the constant function in ${}^{W_+}U_-$ with
value $u$.  It is trivial to check that these functions define Borel
morphisms.  So it suffices to show that there is no Borel morphism
$\xi:\bold W;\bold U\to\bold V$.

Suppose there were such a $\xi$.  Form a forcing extension of the
universe $M$ by first adding a Mathias real $m$ and then, over the
resulting model, adding a Cohen real $c$.  Extend $\xi$ to a Borel
morphism in $M[m,c]$, still called $\xi$, with the same Borel
code.  More precisely, the two functions that constitute $\xi$
correspond (via exponential adjointness in one component) to three
Borel maps
$$
\alpha:V_-\to W_-,\quad
\beta:V_-\times W_+\to U_-,\quad
\gamma:W_+\times U_+\to V_+
$$
with the key property that, for all $x\in V_-$, $w\in W_+$, and $u\in
U_+$, if $W(\alpha(x),w)$ and $U(\beta(x,w),u)$, then
$V(x,\gamma(w,u))$.  We use the same symbols $U_\pm$, $U$, $V_\pm$,
$V$, $W_\pm$, $W$, $\alpha$, $\beta$, and $\gamma$ for the objects in
$M[m,c]$ having the same Borel codes (in $M$).  The key property is a
$\Pi^1_1$ sentence, so it remains true in the forcing extension.

We apply the key property with an arbitrary $x\in V_-\cap M$, with
$w=m$ and with $u=c$.  Then, since $\alpha$ is coded in $M$ and since a
Mathias real dominates all ground model reals, $\alpha(x)\leq^*m$, i.e.,
$W(\alpha(x),m)$.  Similarly, since $\beta$ is coded in the ground
model $M$, the chopped real $\beta(x,m)$ is in $M[m]$ and is therefore
matched by $c$ (since matching $\beta(x,m)$ is a comeager requirement
coded in $M[m]$ and $c$ is a Cohen real over $M[m]$).  That is,
$U(\beta(x,m),c)$.  

By the key property, it follows that $V(x,\gamma(m,c))$.  In other
words, we have a chopped real $\gamma(m,c)\in M[l,c]$ that engulfs
every chopped real $x$ from the ground model $M$.  Recalling the
connection between chopped reals (related by engulfing) and meager
sets (related by inclusion) from Section~3, we find that all the
meager Borel sets coded in $M$ are subsets of a single meager set in
$M[m,c]$.  But Pawlikowski has shown that $M[m,c]$ has no meager set
that includes all the meager Borel sets coded in $M$; see \cite9,
Proposition~1.2 and the discussion following Corollary~1.3.  This
contradiction shows that no such morphism $\xi$ can exist.
\qed\enddemo

The preceding proposition tells us that, of the various combinations
of $\bold U$ and $\bold W$ considered so far, $\bold U;\bold W$ is the
only one admitting a Borel morphism to $\bold V$ and therefore the
only one that can be used to prove
${\bold{add}}(B)\geq\min\{{\bold{cov}}(B),\bb\}$.  In fact,
Corollary~1.3 of \cite9 strongly suggests that $\bold U;\bold W$ is in
some sense equivalent to $\bold V$.  We state this corollary as the
next proposition and give a proof somewhat different from that in
\cite9 (in that we use more forcing but do not use, e.g., the
Kuratowski-Ulam theorem) in order to suggest what the proper sense of
equivalence between $\bold U;\bold W$ and $\bold V$ should be.

\proclaim{Proposition 4 (Pawlikowski)}
Let $M$ be any inner model.  The union of all the meager Borel sets
coded in $M$ is meager if and only if there is a Cohen real $c$ over
$M$ and there is a real $d$ dominating all the reals of $M[c]$.
\endproclaim

\demo{Proof}
The ``if'' half is immediate from the existence of a Borel morphism
$\xi:\bold U;\bold W\to\bold V$ with code in $M$.  Indeed, if $c$ and
$d$ are as in the statement of the proposition, if $x$ is any chopped
real in $M$, and if $\alpha:V_-\to U_-$, $\beta:V_-\times U_+\to W_-$,
and $\gamma:U_+\times W_+\to V_+$ are the parts of $\xi$ (as above),
then we have $U(\alpha(x),c)$ (because $c$ is a Cohen real over $M$
which contains $\alpha(x)$), and $W(\beta(x,c),d)$ (because
$\beta(x,c)$ is in $M[c]$ which $d$ dominates), and therefore
$V(x,\gamma(c,d))$ (because $\xi$ is a morphism).  Thus, $\gamma(c,d)$
engulfs all chopped reals $x\in M$, and therefore it codes a meager
set that includes all meager Borel sets coded in $M$.

For the ``only if'' direction, suppose we have a meager set that
includes all the meager Borel sets coded in $M$.  Then, by Section~3,
we have a chopped real $(y,\Theta)$ that engulfs all the chopped reals
of $M$.  Then clearly $y$ matches every chopped real from $M$ and is
therefore Cohen generic over $M$.  It remains to produce a real $g$
dominating all the reals of $M[y]$.  We claim that such a $g$ is
obtained by letting $g(n)$ be the right endpoint of the next interval
in $\Theta$ after the one that contains $n$.

To verify that this $g$ works, we consider any $f:\omega\to\omega$ in
$M[y]$ and show that $f\leq^*g$.  We may suppose without loss of
generality that $f$ is non-decreasing.  Since $f$ belongs to the Cohen
extension $M[y]$, it is the denotation with respect to $y$ of some
name $\dot f\in M$.  (Here and in the rest of this proof, names and
forcing are with respect to the usual notion of forcing for adding a
Cohen real, ${}^{<\omega}2$ ordered by reverse inclusion.)  We assume
without loss of generality that all conditions force ``$\dot f$ is a
function from $\omega$ to $\omega$.''
 
We define a chopped real $(x,\Pi)\in M$ with the property that, if
$[a,b]$ is any interval in $\Pi$ and if $p$ is any Cohen condition
that (has length at least $b$ and) agrees with $x$ on $[a,b]$, then
$p$ forces ``$\dot f(a)\leq b$.''  We proceed by induction.  After $n$
intervals of $\Pi$ and the restrictions of $x$ to those intervals have
been defined, we produce the next interval $[a,b]$ and the restriction
of $x$ to it as follows.  Of course $a$ is the first number not in the
intervals already defined.  Fix a list $u_0,u_1,\dots,u_{r-1}$ of all
the functions $a\to2$.  We inductively define functions $x_i:[a,l_i)$
with $a=l_0\leq l_1\leq\dots\leq l_r$ by starting with the empty
function as $x_0$ and obtaining $x_{i+1}$ as an extension of $x_i$
such that $u_i\cup x_{i+1}$ forces a particular value $v_i$ for $\dot
f(a)$.  Such an extension exists because $\dot f$ is forced to be a
total function on $\omega$.  After $x_r$ and $l_r$ have been reached,
let $b$ be the largest of $l_r$ and all the $v_i$, and extend $x_r$
arbitrarily to $[a,b]$; this defines the restriction of the desired
$x$ to the next interval $[a,b]$ of $\Pi$.  If a condition $p$ agrees
with $x$ on $[a,b]$ then it agrees with some $u_i$ on $a$ and
therefore extends $u_i\cup x_{i+1}$ and forces ``$\dot f(a)=v_i\leq b$''
as required.  

As $(x,\Pi)$ is in $M$, it is engulfed by $(y,\Theta)$.  That is, each
interval $[m,n]$ of $\Theta$, with only finitely many exceptions,
includes an interval $[a,b]$ of $\Pi$ on which $x$ and $y$ agree.
Then the initial segment $y\res b$ of $y$ is a condition of the sort
considered in defining $x$, so it forces ``$\dot f(a)\leq b$.''  Since $y$
is Cohen generic, this forced statement is true in $M[y]$, i.e.,
$f(a)\leq b$.  But then, as $f$ is non-decreasing, we have, for all
elements $k$ of the interval of $\Theta$ immediately preceding
$[m,n]$, that 
$$
f(k)\leq f(a)\leq b\leq n=g(k).
$$
This applies to all sufficiently large $k\in\omega$ because $[m,n]$
can be any interval of $\Theta$ with only finitely many exceptions.
Therefore, $f\leq^*g$, as required.
\qed\enddemo

We would like to regard the ``only if'' half of the preceding proof as
presenting a morphism $\eta:\bold V\to\bold U;\bold W$ in the
direction opposite to the $\xi$ involved in the ``if'' half.  This
would mean that the proof involves (1) a contruction $\eta_-$ whose
input consists of a chopped real $(z,\Psi)$ (from $U_-$) and a
function $\Phi:U_+\to W_-$ and whose output is a chopped real
$(x,\Pi)$ and (2) a construction $\eta_+$ whose input consists of a
chopped real $(y,\Theta)$ and whose output consists of a real $c$ and
an element $g$ of $W_+$.  The required key property for a morphism is
that (with the notation of (1) and (2)), if $(y,\Theta)$ engulfs
$(x,\Pi)$, then $c$ matches $(z,\Psi)$ and $g$ dominates $\Phi(c)$.

Our proof looks vaguely but not exactly like this.  Notice that a
Cohen-forcing name $\dot f$ of a real, as in our proof of the
proposition, defines a function $\Phi$ into ${}^\omega\omega$ from a comeager
subset $D$ of ${}^\omega2$, namely, $\Phi(u)(n)$ is the unique number
forced to be the value of $\dot f(n)$ by some initial segment of $u$.
($D$ consists of those $u$ whose initial segments force values for all
$\dot f(n)$; this $D$ is comeager because $\dot f$ is forced to name a
real.)  For a $\Phi$ of this special sort (or, more precisely, for any
extension of such a $\Phi$ to all of ${}^\omega2=U_+$), our proof
produced, from $\dot f$, a certain $(x,\Pi)$.  For the output of
$\eta_-$, when the input consists of $(z,\Psi)$ and such a special
$\Phi$, we take a chopped real that engulfs $(z,\Psi)$, the $(x,\Pi)$
constructed from $\dot f$ in the proof, and an $(x',\Pi')$ such that
$\m x',{\Pi'}$ is included in the proper domain $D$ of $\Phi$.  (It is
trivial to produce a chopped real engulfing any finitely (or even
countably) many given chopped reals.)  For the
output of $\eta_+$ on input $(y,\Theta)$, we take the pair whose first
component is $y$ and whose second component is the $g$ constructed
from $\Theta$ during the proof (the ``end of the next interval''
function).  Then the key property that $\eta$ needs is just what is
established by the proof.

All this, however, was done only for the very special $\Phi$'s that
correspond to names $\dot f$.  The proof gives no hint what $\eta_-$
should do if its second argument $\Phi$ is not of this form.  This
difficulty is, of course, connected with the difficulty mentioned
earlier that sets like ${}^{U_+}W_-$ are one type higher than the
(Borel) sets of reals that our theory is equipped to handle.  

Both difficulties can be attacked by working not with all functions
$U_+\to W_-$ but with some subfamily that can be coded by reals (or
with the codes rather than the functions).  In the case at hand, it is
tempting to take this family to be just those $\Phi$'s that are given
by Cohen names $\dot f$; these names, which are essentially reals, can
then serve as the codes.  But this choice of a subfamily is obviously
tailored to just this one example.  We can do much better by noticing
that any Borel function $\Phi:U_+\to W_-$ agrees on some comeager set
$A$ with the function given by some $\dot f$.  Thus, we can allow
Borel functions as inputs to $\eta_-$ at the cost of redefining the
output of $\eta_-$ so that it also engulfs a chopped real whose
matching set is included in $A$.

The restriction to Borel functions is one that can be sensibly imposed
in general.  That is, we can modify the general definition of
sequential composition $\bold A;\bold B$, when all the components of
$\bold A$ and $\bold B$ are Borel sets and relations, by replacing
the set ${}^{A_+}B_-$ of functions with the subset of Borel functions
or, better, with the set of Borel codes for such functions.

To be specific about the coding, we note that Borel functions are
precisely the functions recursive in a real and the type 2 object
${}^2E$ (see for example \cite6, Theorem~VI.1.8).  So as a code for
such a function we can use the pair of the Kleene index (a natural
number) for this recursive function and the real parameter.

This framework covers the standard examples, including the particular
ones we have discussed.  It should, however, probably be extended a
bit because the set of codes of Borel maps from a Borel set to a Borel
set is in general not itself a Borel set; it is only
$\boldsymbol\Pi^1_1$.  So it seems reasonable to allow the domains and
codomains of relations to be $\boldsymbol\Pi^1_1$ (i.e.,
semi-recursive in ${}^2E$ and a real) sets of reals; the relations
should probably still be required to be Borel in the weak sense that
they and their complements relative to the domain and codomain are
$\boldsymbol\Pi^1_1$ (so that dualization works).  Morphisms should be
functions with $\boldsymbol\Pi^1_1$ graphs, i.e., partial recursive in
${}^2E$ and reals but total on the $\boldsymbol\Pi^1_1$ sets in
question.

Rather than speculate further on the basis of very limited examples,
let me just list what I would like to see in an ideal framework.  The
sets, relations, and morphisms should come from a class that is broad
enough to encompass the known examples but narrow enough to have
absoluteness properties so that non-existence of incorrect morphisms
can be proved.  Furthermore, the framework should be closed under the
naturally occurring constructions, including a suitable version of
sequential composition.  I believe that the framework described in the
preceding paragraph may do all this, but this has not yet been fully
checked.

I should perhaps mention that Pawlikowski told me that, with suitable
coding, the components of morphisms that occur in classical proofs of
cardinal characteristic inequalities can be taken to be not merely
Borel but continuous.  It is not clear, however, that this helps in
the present situation, since the (natural) codes for continuous
functions form only a $\boldsymbol\Pi^1_1$ set.

\head
6. Unsplitting
\endhead

In this section, we briefly discuss a situation where sequential
composition arises naturally in connection with an inequality
involving just two cardinal characteristics, and a rather trivial
inequality at that.  Recall that the unsplitting number $\rr$ is the
smallest number of infinite subsets of $\omega$ such that no single
set splits them all into two infinite pieces.  Equivalently, it is the
norm of $\bold R=({}^\omega2,\Cal P_\infty(\omega),R)$ where $R(f,Y)$
means that $f$ is almost constant on $Y$, i.e., is constant on $Y-F$
for some finite $F$.  

We define a similar characteristic $\rr_3$ using splittings into three
pieces rather than two.  That is, $\rr_3$ is the norm of $\bold
R_3=({}^\omega3,\Cal P_\infty(\omega),R_3)$ where again $R_3(f,Y)$
means that $f$ is almost constant on $Y$.

It is easy to see that $\rr_3=\rr$.  In fact, there are morphisms
$\xi:\bold R_3\to\bold R$ and $\eta:\bold R;\bold R\to\bold R_3$
defined as follows.  $\xi_-$ is the inclusion map
${}^\omega2\to{}^\omega3$ and $\xi_+$ is the identity map of $\Cal
P_\infty(\omega)$.  To define $\eta$, we first introduce, for
$f:\omega\to3$, the notations $f'$ and $f''$ for the functions
$\omega\to2$ obtained by identifying 2, whenever it occurs as a value
of $f$, with 0 or 1 respectively.  That is,
$$ 
f':\omega\to2:n\mapsto
\cases
0,&\text{if $f(n)=0$}\\
1,&\text{otherwise,}
\endcases
$$
and $f''$ is defined similarly except that the first case is ``if
$f(n)=0$ or 2.''  Also fix, for each infinite $Y\subseteq \omega$ a
bijection $e_Y:\omega\to Y$, for example the unique increasing
bijection.  Now we define $\eta$ as follows.  $\eta_-$ sends
$f:\omega\to3$ to the pair consisting of $f'$ and the map $G:\Cal
P_\infty(\omega)\to{}^\omega2$ defined by $G(Y)=f''\circ e_Y$.  (In
other words, $G(Y)$ is essentially $f''\res Y$ but with its domain
shifted from $Y$ to $\omega$ by $e_Y$.)  $\eta_+$ sends a pair of
infinite sets $(Y,Y')$ to $e_Y(Y')$, the set that ``occupies in $Y$
the locations that $Y'$ occupies in $\omega$.''  It is easy to check
that $\xi$ and $\eta$ are morphisms, and therefore
$\rr\leq\rr_3\leq\rr\cdot\rr=\rr$.

The idea behind the morphism $\eta$ is that to get a 3-unsplit family
(in the obvious sense), it suffices to start with a 2-unsplit family
$\Cal R$ and then, within each of its members $\hskip0.01cm Y$, form a
new 2-unsplit family by transferring $\Cal R$ from $\omega$ to $Y$ via
$e_Y$.  The union of these new, transferred families is 3-unsplit
because, if $f:\omega\to3$, then there is some $Y\in\Cal R$ on which
$f$ takes at most two values infinitely often and there is some $Y'\in
e_Y[\Cal R]$ on which those two values are reduced to one.

Notice that the ``second order'' families $e_Y[\Cal R]$ depend in an
essential way on the $Y$'s.  In other words, the domain of $\eta$ in
this proof apparently needs to be a sequential composition.  There is
no evident way to use a product or even an old product instead.

We therefore conjecture that there is no Borel morphism from the old
product of two copies of $\bold R$ to $\bold R_3$.  Intuitively, this
means that, in the two uses of a 2-unsplit family $\Cal R$ to produce
a 3-unsplit family, the second use cannot be made independent of the
first.

Although the conjecture seems highly plausible, we do not know how to
prove even the weaker conjecture that there is no Borel morphism from
$\bold R$ to $\bold R_3$.  The intuitive meaning of this is merely
that, to produce a 3-unsplit family, we must use the 2-unsplit family
$\Cal R$ twice, not just once.

To see the difficulty in proving these conjectures, recall that our 
previous proofs of non-existence of Borel morphisms involved the
construction of suitable forcing models.  To apply this method
directly to prove the conjecture, we would want to find a forcing
extension that contains an infinite subset $Y$ of $\omega$ unsplit by
all ground model functions $\omega\to2$ but split by some ground model
function $f:\omega\to3$.  This clearly cannot be achieved, since a set
unsplit by $f'$ and $f''$ (in the notation of the definition of $\eta$
above) is also unsplit by $f$.

An alternative approach would be concentrate on the dual morphism
$\bold R_3{^\bot}\to\bold R^\bot$.  Now the forcing method suggests
building a forcing extension in which every function $\omega\to2$ is
almost constant on an infinite $Y$ from the ground model but some
$f:\omega\to3$ is not.  But this cannot be achieved either.  If
$f:\omega\to3$, then there is an infinite $Y$ in the ground model on
which $f'$ is almost constant.  Furthermore, $f''\circ e_Y$ is almost
constant on some $Y'$ in the ground model.  If $e_Y$ was
defined reasonably (e.g., as the increasing bijection), then it too
is in the ground model, and so is $e_Y(Y')$ on which $f$ is almost
constant.  

So the forcing method seems to be of no use in proving even the weak
form of the conjecture.  Some direct analysis of Borel morphisms seems
to be needed.

\enddocument